\documentclass[]{amsart}
\usepackage{amssymb, amsmath, amsthm, amscd,ifthen}
\usepackage[dvips]{graphics}
\usepackage{graphicx}
\usepackage{epsfig}

\begin{document}

\hfill{MSC 52A20, 52A38, 52B11, 52B55, 49Q20, 65D17, 68U07}

\bigskip

\centerline{\bf Blaschke Addition and Convex Polyhedra}

\centerline{\bf Victor Alexandrov\footnote{Partially supported by
the Russian Foundation for Basic Research (grant 03--01--00104) and the Ministry
for Education of the Russian Federation (Grant E02--1.0--43).},
Natalia Kopteva, S.~S.~Kutateladze}

\bigskip

Abstract. This is an extended version of a talk on
October 4, 2004
at the research seminar ``Differential geometry and applications''
(headed by Academician A.~T. Fomenko) at  Moscow State University.
The paper contains an overview of available (but far from well-known)
results about the Blaschke addition of convex bodies, some new theorems
on the monotonicity of the volume of convex bodies (in particular,
convex polyhedra with
parallel faces) as well as description of a software for visualization
of polyhedra
with prescribed outward normals and  face areas.

\bigskip

{\bf 1. The vector area of the surface of a polyhedron.}
In this article we take the following fact as well-known:
{\sl Let $P$ be a compact convex polyhedron in $\Bbb R^m$ ($m\geq 2$),
$n_1,\dots,n_k$ be unit vectors of outward normals to its $(m-1)$-dimensional
faces and let
$F_1,\dots,F_k$ be $(m-1)$-dimensional volumes of its $(m-1)$-dimensional
faces. Then
$$
\sum\limits_{j=1}^{k}F_jn_j=0.\eqno(1)
$$}

Sometimes this fact is stated as follows:
``The vector area of the surface of a convex polyhedron equals zero''
(see, for example, \cite{AD50}).
From a hydrodynamical point of view this means that a polyhedron
submersed in a fluid  under only the action of the forces
of the fluid pressure
remains in equilibrium.
Convexity of the polyhedron is of course a redundant requirement.
Item 1 is a direct consequence of the Stokes theorem.
Its ``smooth'' analog says that the integral of the unit outward normal
over a closed surface equals zero. This t is given, for example,
in a celebrated calculus problem book by B.~P.~Demidovich (problem 4381).
We do not prove Item 1 since it is obvious,
but we give an analog of it for a domain on a sphere \cite{Al02}:

{\sl Let $\Bbb S^2$ be a unit sphere in $\Bbb R^3$ and
let $D$ be a domain in $\Bbb S^2$ with piecewise smooth boundary.
Let $\hat N$ stand for a vector field on the sphere $\Bbb S^2$
that maps each point of the sphere to a unit outward normal
with  base this point of the sphere .
Let $\hat n$ denote a vector field on the smooth part of the boundary
$\partial D$ of $D$ that maps each point to a unit outward
(with respect to $D$)
vector which is normal to $\partial D$ and parallel to the plane
tangent to the sphere.
Let $\sigma$ denote the standard measure on $\Bbb S^2$ and
$s$ denote the arc length on $\partial D$. Then
$$
2\iint_D \hat N\, d\sigma +\int_{\partial D} \hat n\, ds =0.\eqno(2)
$$}

It is natural to think that the first integral in (2)
gives the barycenter of the domain $D$ on the sphere.
It would be interesting to find, by using  (2),
the ``geographical centers'' of
domains on the Earth whose curvature cannot be neglected
(for example, Russia).

Another basic point of this survey is the following celebrated
theorem of H.~Minkowski \cite{AD50}.

{\bf 2. Theorem (Minkowski)}.
{\sl Let $m\geq 2$, $F_1,\dots, F_k$ be positive real,
$n_1,\dots,n_k$ be unit vectors in $\Bbb R^m$ not lying
in a single hyperplane
such that $\sum_{j=1}^{k}F_jn_j=0.$
Then there exists a compact convex polyhedron $P$ in $\Bbb R^m$ such
that the vectors
$n_1,\dots,n_k$ (and only they) are the outward normals to
$(m-1)$-dimensional
faces of $P$ and $F_1,\dots, F_k$ are the $(m-1)$-dimensional volumes
of the $(m-1)$-dimensional faces of $P$.
Moreover, such a polyhedron $P$ is unique up to translation.}

In \cite{Al04}, this theorem is abstracted to the class of
nonconvex polyhedra
in $3$-dimensional Euclidean space which are called herissons.
The definition of  herisson  in \cite{Al04} is rather complicated
and requires
injectivity of the spherical map of a nonconvex polyhedron.
Referring to a recent talk of N.~P.~Dolbilin, A.~T.~Fomenko suggested
to significantly simplify the definition of  herisson
as follows.
A herisson is defined to be a set of vectors not lying in
a single hyperplane
and equipped with some numbers  meeting (1).

We can define  addition on the set of such herissons
in a natural way.
The sum of two herissons $\{n_1,\dots,n_k;F_1,\dots, F_k\}$ and
$\{n'_1,\dots,n'_l;F'_1,\dots, F'_l\}$
is defined to be the set of vectors
${\mathcal N}=\{n_1,\dots,n_k\}\cup\{n'_1,\dots,n'_l\}$ together
with the positive numbers determined as follows:

$\bullet$ if a vector $n\in\mathcal N$ belongs to only one of the sets
$\{n_1,\dots,n_k\}$ or $\{n'_1,\dots,n'_l\}$ and equals, say, $n_j$,
then assign the number $F_j$ to the vector $n$;

$\bullet$ if $n\in\mathcal N$ belongs to both
$\{n_1,\dots,n_k\}$ and $\{n'_1,\dots,n'_l\}$ and
equals, say, $n_j$ and $n'_p$,
then assign the number $F_j+F'_p$ to the vector $n$.

It is clear that the set of vectors $\mathcal N$ does not
lie in a hyperplane and
the condition (1) holds, that is the sum of herissons is a herisson.

By the fact 1, for each compact convex polyhedron, there exists a herisson
whose vectors are the unit outward normals to the faces of the highest
dimension equipped with ``areas'' of the faces.
Two polyhedra give rise to two herissons. The sum of these herissons
is a herisson,
which corresponds to a compact convex polyhedron by the Minkowski theorem.
Such a polyhedron is determined up to translation.
So we can define a new addition on the set of classes of translates of convex polyhedra.
Examples below show that this operation is differed from
the Minkowski sum of polyhedra.

Considering the importance of the Minkowski sum of convex bodies,
it is natural to study the properties of this new operation.
In order to complete the picture,
let us first recall some basic facts about the Minkowski sum.

{\bf 3. Minkowski sum of convex bodies} \cite{AD50,BF02,Le85,Ly56,Sc93a}.
The set $K+L=\{ z\in\Bbb R^m | z=x+y, \ x\in K,\ y\in L\}$
is called the Minkowski (or vector) sum  of nonempty convex compact sets
$K,L\subset\Bbb R^m$.

Another definition of  Minkowski sum uses the support function
of a convex set $K$
which is the function $h_K$ defined on $\Bbb R^m$ by the formula
$$
h_K(x)=\sup\limits_{y\in K} (x,y),
$$
where $(x,y)$ stands for the usual inner product of vectors $x$ and $y$
in $\Bbb R^m$.

The support function of a convex compact set is a positive homogeneous
convex function.
Moreover, there is a one-to-one correspondence between
positive homogeneous convex functions
and convex compact sets in $\Bbb R^m$.
Let us give another definition of the Minkowski sum
(clearly equivalent to the previous )
which is based on the bijectivity of this correspondence:
the Minkowski (or vector) sum of nonempty convex compact sets $K$ and $L$ from $\Bbb R^m$
is the convex set $K+L$ whose support function $h_{K+L}$ is
the sum of the support functions
of $K$ and $L$; i.e., $h_{K+L}=h_K+h_L$.

It is not difficult to show that the sum of some translates (i.e.,
images of a convex
body under some translation of the space) $K$ and $L$ is a translate
of~$K+L$. Therefore, we can speak correctly about the Minkowski sum of
translates of convex bodies.
Allowing some ambiguity, we will sometimes identify the class of
translates of a
convex body with the body itself.

The notion of the vector sum of convex bodies was introduced by Minkowski
in study of some questions about the isoperimetric inequality \cite{Mi03}.
Since then, this notion was well studied and found a wide application to
various areas of mathematics; for example, see \cite{Ga02}
and \cite{Sc93a}.
The most known result  using the Minkowski sum is the Brunn--Minkowski
inequality
$$
(\mbox{Vol\, }(K+L))^{1/m}\geq (\mbox{Vol\, }K)^{1/m}+(\mbox{Vol\, }L)^{1/m},\eqno(3)
$$
which holds for all convex compact sets $K,L\subset\Bbb R^m$.

{\bf 4. Blaschke sum of convex bodies.}
Let us come back to the addition of convex polyhedra which was
described in Section 2.
The terminology is questionable to some degree in fact.
For example, B. Gr\"unbaum says \cite[p. 339]{Gr03}:
``the process we have called
the Blaschke addition was first described by Blaschke \cite{Bl67}
for smooth convex sets, although, even earlier, the
corresponding addition of polytopes occurs implicitly in the work
of Minkowski \cite{Mi87}.''

The Blaschke sum is defined for arbitrary convex bodies in $\Bbb R^m$;
i.e.,  compact convex sets with nonempty interior.
(In fact, we may omit the claim of nonempty interior, but this
leads to some technical difficulties that we are reluctant
to discuss in this paper.)
To this end ,we need the following definition:
the surface area measure $S(K,\cdot)$ of a convex body $K$ in $\Bbb R^m$
is the additive set function on the unit sphere $\Bbb S^{m-1}$ which
associates to each
subset $U\subset\Bbb S^{m-1}$ the (surface) measure of
the set of the points $x\in\partial K$
admitting a unit outward normal to $\partial K$
such that
if we shift this vector to make its base coincide with the center of
$\Bbb S^{m-1}$ then its endpoint will belong to $U$.

Let us give some standard examples of surface area measures.

(4.1) The surface area measure of a sphere of radius $r$
in $\Bbb R^m$ is proportional to the surface area measure
of the unit sphere
$\Bbb S^{m-1}$ with rate $r^{m-1}$.

(4.2) The surface area measure of a polyhedron $\Bbb R^m$
is concentrated
in the normals to the faces.
Its value at the ``endpoint'' $x\in\Bbb S^{m-1}$
of such a normal equals the $(m-1)$-dimensional volume of the face.

(4.3) The surface area measure of a $C^2$-smooth convex body
is given by the
identity
$$
S(K,U)=\int_U \frac{d\sigma}{\varkappa},
$$
where $\varkappa$ is the Gaussian curvature of the surface $\partial K$
and $\sigma$ is the standard measure on the sphere $\Bbb S^{m-1}$.

The notion of  surface area measure of a convex body was introduced in
1938 by A.~D. Alexandrov \cite{AD38} as well as by W.~Fenchel and
B.~Jessen \cite{FJ38}
independently.
This notion is a typical ingredient of  convex analysis
and allows one to  work uniformly
with smooth and nonsmooth convex bodies.

We we are ready  now to give the central definition of this article.

The {\it Blaschke sum} (of the classes of translates) of two convex bodies
$K$ and $L$ in $\Bbb R^m$ is
the (class of translates of the) convex body $K\# L$ whose
surface area measure is
the sum of surface area measures of the summands:
$S(K\# L,\cdot)=S(K,\cdot)+S(L,\cdot)$.

Correctness of the definition of  Blaschke sum rests on the following
form of the Minkowski theorem on unique existence of a convex body
with a given surface area measure which was first proved by 
A.~D. Alexandrov \cite{AD38}:

{\bf 4.4. Theorem.} {\sl Let $\mu$ be an additive set function on the unit sphere
$\Bbb S^{m-1}\subset\Bbb R^m$ such that
$$
\int_{\Bbb S^{m-1}}x\, d\mu=0
$$
(here $x$ denotes the radius-vector of a point on $\Bbb S^{m-1}$)
and there is no hyperplane $\Pi$ such that
$\mu$ is concentrated in $\Pi\cap\Bbb S^{m-1}$.
Then there exists a convex body $K\subset\Bbb R^m$
 whose surface area measure coincides
with $\mu$; i.e,, $S(K,\cdot)=\mu(\cdot)$.
Moreover, this body $K$ is unique up to translation.}

Theorem 2 is a particular case of Theorem 4.4.

The Blaschke sum of polyhedra coincides with
the addition of herissons as described in Section 2.

{\bf 5. Properties of the Blaschke sum.}
Although an analogy between the definitions of  Minkowski and
Blaschke sums is
obvious, the Blaschke sum has been studied to a much lesser degree.
Let us discuss the properties of the Blaschke sum which we
managed to find in the literature.

(5.1) {\sl On the Euclidean plane $\Bbb R^2$, the Blaschke sum of two
convex polygons
coincides with their Minkowski sum: $K\# L=K+L$.}

This can be proved  for instance by inspecting a detailed construction
of the Minkowski
sum of convex polygons in \cite{Ly56}.

(5.2) {\sl In dimensions greater than or equal to $3$,
the Blaschke sum of convex polyhedra
does not  coincide necessarily with their Minkowski sum.}

Skipping details,
we give the following picture from \cite{AD50} as an example
exhibiting the Minkowski sum of two regular tetrahedra one of which is
obtained from the
other by rotation through $90^\circ$ around the vertical axis.
From Figure~1 it is clear that the Minkowski sum of these tetrahedra
is a polyhedron with
14 faces, while their Blaschke sum is a polyhedron with 8 faces by
definition.
The difference transpires.

\begin{figure}
\centering
\includegraphics[width=8 cm]{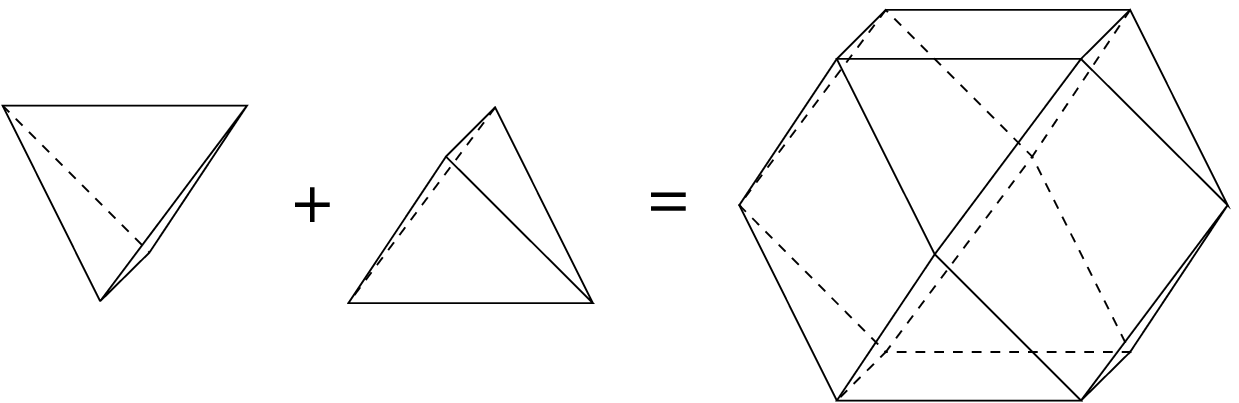}
\caption{The Minkowski sum of two tetrahedra}\label{minkowski}
\end{figure}

(5.3) {\sl For all convex bodies $K,L$ in $\Bbb R^m$, the
Kneser--S\"uss inequality holds:
$$
(\mbox{Vol\, }(K\#L))^{1-1/m}\geq (\mbox{Vol\, }K)^{1-1/m}+(\mbox{Vol\, }L)^{1-1/m},
$$
moreover,  equality takes place if and only if $K$ and $L$ are
homothetic.}

A proof can be found in \cite[Theorem 7.1.3]{Sc93a} and \cite{KS32}.

(5.4) {\sl Every compact convex polyhedron with nonempty interior in
$\Bbb R^m$
can be expressed as the Blaschke sum of finitely many simplices
in $\Bbb R^m$.
Furthermore, these simplices can be chosen so that
the sum contains at most
$f-m$ summands, where $f$ is the number of faces of the polyhedron.}

A proof can be found in \cite[pp. 334--335]{Gr03} and \cite{FG64}.

The Blaschke addition is not adequately explored. This is
confirmed
for example by the fact that
it is still unknown whether
we can generate functionals analogous to
the mixed volumes generated by the Minkowski sum
by using the Blaschke sum
and  prove inequalities for these functionals by analogy to
the Minkowski and Alexandrov--Fenchel inequalities for mixed volumes.

\begin{figure}
\begin{minipage}[b]{0.50\linewidth}
\centering
\includegraphics[width=3.2 cm]{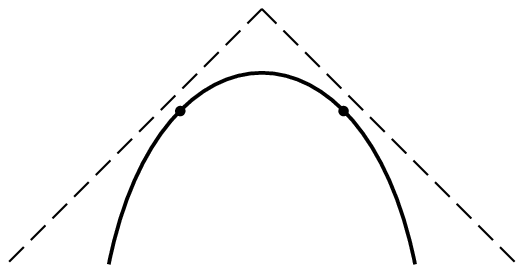}\\
\bigskip
a)
\end{minipage}\hfill
\begin{minipage}[b]{0.50\linewidth}
\centering
\includegraphics[width=3.2 cm]{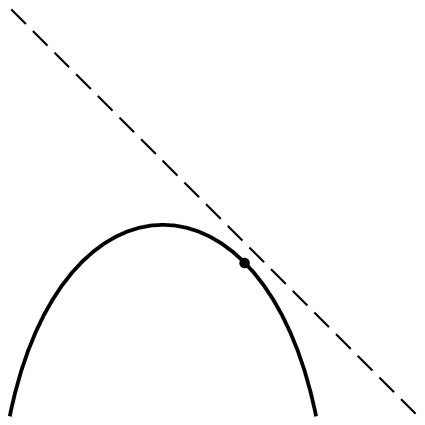}\\
\bigskip
b)
\end{minipage}\hfill
\caption{Maximization of a concave functional}\label{minimum}
\end{figure}

{\bf 6. The Blaschke sum and extremal problems.}
It is known that the classical isoperimetric problem consists in
to finding
a body of maximal volume
among all convex sets in $\Bbb R^m$ of a given surface area.
In order to apply the tools of analysis,
the problem must be  parametrized for example
by identifying convex sets with their support functions.
This leads to the problem of finding a maximum of a concave
functional, namely, the volume of the set, in a vector space of
continuous functions
on the sphere, namely, support functions, subject to a constraint,
namely, a fixed surface area.
Difficulties in the proof of existence and uniqueness of a solution
of the classical isoperimetric problem appear because
both the constraint
( surface area) and the  objective functional  ( volume) are concave
(see Figure~\ref{minimum}~a)).
Concavity of both functionals is guaranteed
exactly by the Brunn--Minkowski inequality (3)
(there is a corresponding analog for
the surface area).

These difficulties naturally disappear if the problem is parametrized by
identifying  convex bodies with their surface area measures.
Then the constraint (that the surface area be constant) turns into
a hyperplane in the vector space of all surface area measures and the
objective
functional ( volume) remains concave which is guaranteed by the
Kneser--S\"uss inequality (5.3).
In this situation, we can use the conventional methods of
the theory of convex programming
so that both existence and uniqueness of a solution is obtained
``for free'' (see Figure~\ref{minimum}~b)) and a solution can be found
from the Euler equation which says that in the point of extremum
the tangent  plane to the objective functional
is parallel to the tangent plane to the constraint.

In \cite{Ku73}, it was studied why the established modern theory of
extremal problems helps little in solving isoperimetric problems.
In that paper, a functional-analytical view of isoperimetric problems
was developed, in particular, it was suggested to subdivide
the problems into those that are
``strategically right'' to be parametrized by support functions
(then the Minkowski sum plays a key r\^ole), by surface area measures
(then the Blaschke sum is of fundamental importance), and those
that do not admit any parametrization that leads to convex extremal
problems.
Let us explain this by examples of some  problems.

(6.1) {\sl The Urysohn problem}:
Find a body of maximal volume
among convex bodies in $\Bbb R^m$ of diameter at most 1.

From a functional-analytical point of view, it is a question of
maximization of a concave functional, namely, the volume,
on a convex set in the vector space of support functions
(indeed, the fact that the width of a convex body $K$ in direction
$n\in\Bbb R^m$ is at most 1 amounts to
the linear inequality $h_K(n)-h_K(-n)\leq 1$).
By a general argument, the problem has a unique solution
that can be found
from the Euler equation. A solution to the problem was found by
P.~S.~Urysohn in 1924~\cite{Ur24}. He  proved that the sought body is
exactly a ball.

(6.2) {\sl Inner isoperimetric problem}:
Find a convex body of maximal volume among convex bodies in
$\Bbb R^m$ of a prescribed surface area
which lie in a given convex body $T$.

From a functional-analytical point of view,  complexity of this problem
consists in the fact that the set of convex surfaces
lying a given convex body
is convex with respect to the Minkowski addition (but not with respect
to the Blaschke addition), while the surface area is linear with respect
to the Blaschke addition (and with respect to the Minkowski addition,
it is just concave in much the same way as  volume).
Such a ``mixture of styles'' when one part of conditions behaves well
with respect to the Minkowski sum and another part of conditions
behaves well
with respect to the Blaschke sum essentially complicates the problem
and results in the fact that the inner isoperimetric problem is
less studied
even in the case when $T$ is a tetrahedron.
From the necessary extremum condition (the Euler equation),
it is clear that the boundary of the extremal body consists of
flat pieces (where the body contacts  the boundary of the tetrahedron)
and pieces of surfaces of constant mean curvature
(where its surface is inside the tetrahedron).
In 1994, A.~V. Pogorelov showed that for every tetrahedron $T$,
not necessarily regular,
and every positive number $H>1/R$, where $R$ is the radius of the ball
inscribed in the tetrahedron $T$,
there exists a smooth closed convex surface inside $T$ that contacts
with all of its faces such that the part of the surface contained
strictly within $T$
has a constant mean curvature $H$ \cite{Po94}.
This surface can be described in the following way:
take the Blaschke sum of the tetrahedron $T$ with a ball of a special
radius and then take the Minkowski sum of the Blaschke addition result
with another ball \cite{Po94}.
However, a surface constructed in this way satisfies only the
necessary extremum condition.
It is still unclear whether it is indeed a solution to the inner
isoperimetric
problem and if the problem has a unique solution.
There are no other significant advances in the solving of the
inner isoperimetric
problem.

The Blaschke sum was systematically used by S.~S. Kutateladze
\cite{Ku73,Ku98,Ku02}
for solving  isoperimetric problems with many subsidiary constraints.

We mention also that various generalizations of the Blaschke sum of
convex compacts
were considered by W.~Firey \cite{Fi61,Fi65} and other authors
\cite{Fe78,Fe79,Sc93b}.

{\bf 7. New theorems on monotonicity of the volume of convex bodies.}

{\bf 7.1. Theorem.} {\sl Let $K$ and $L$ be convex bodies in $\Bbb R^m$ ($m\geq 2$)
and let the surface area measure of $K$ do not exceed
the surface area measure of
$L$, that is $S(K,\cdot)\leq S(L,\cdot)$.
Then the volume of $K$ does not exceed the volume of $L$, that is
$\mbox{Vol\, } K\leq\mbox{Vol\, }L$.}

{\bf Proof.}
Take $t\in(0,1)$. Consider an additive set function $\mu(\cdot)$
of the unit sphere $\Bbb S^{m-1}$ such that
$\mu(\cdot)=S(L,\cdot)-t S(K,\cdot)$.
Since $t<1$, for every subset $U\subset\Bbb S^{m-1}$ such that
$S(K,U)>0$, we have $\mu(U)>0$.
Since $S(K,\cdot)$ is not concentrated in a hyperplane
(more precisely, its support is not contained in $\Pi\cap\Bbb S^{m-1}$
for all hyperplanes $\Pi$),
we conclude that $\mu$ is not concentrated in a hyperplane.
Finally, the condition
$$
\int_{\Bbb S^{m-1}} x\, d\mu=0
$$
holds, because analogous integrals with respect to additive set functions
$S(K,\cdot)$ and $S(L,\cdot)$ vanish.

By the Minkowski theorem 4.4, there exists a convex body $M$
(unique up to translation)
in $\Bbb R^m$ whose surface area function coincides with $\mu$,
that is $\mu(\cdot)=S(M,\cdot)$.
But then $S(L,\cdot)=S(M,\cdot)+ tS(K,\cdot)$, and so
$L=M\#(t^{1/(m-1)}K)$.
By the Kneser--S\"uss inequality, we then have
$$
(\mbox{Vol\, }L)^{1-1/m} = (\mbox{Vol\, }(M\#(t^{1/(m-1)}K)))^{1-1/m}\geq
$$
$$
\geq (\mbox{Vol\, }M)^{1-1/m}+ (\mbox{Vol\, }(t^{1/(m-1)}K))^{1-1/m}=
$$
$$
=(\mbox{Vol\, }M)^{1/(m-1)}+t\, (\mbox{Vol\, }K)^{1-1/m}\geq t\, (\mbox{Vol\, }K)^{1-1/m}.
$$
Tending $t$ to $1$, we get $(\mbox{Vol\, }L)^{1-1/m}\geq (\mbox{Vol\, }K)^{1-1/m}$,
which completes the proof.

Let us give some obvious corollaries from Theorem~7.1.

{\bf 7.2. Corollary.} {\sl Let $\{n_1,\dots,n_k\}$ be the unit
outward normals to
$(m-1)$-dimensional faces of a compact convex polyhedron
$P\subset\Bbb R^m$ with nonempty interior and let
$\{n'_1,\dots,n'_l\}$ be the unit outward normals to
$(m-1)$-dimensional faces of a compact convex polyhedron
$Q\subset\Bbb R^m$ with non-empty interior.
Put ${\mathcal N}=\{n_1,\dots,n_k\}\cup\{n'_1,\dots,n'_l\}$.
Given $n\in\mathcal N$, let $\Gamma_P(n)$ be the face
of $P$ with outward normal $n$
and let $\Gamma_Q(n)$ be the face of $Q$ with outward normal $n$.
At least one of these faces is of dimension $m-1$.
Suppose that for every vector $n$ in $\mathcal N$,
the $(m-1)$-dimensional volume of $\Gamma_P(n)$
is not less than
the $(m-1)$-dimensional volume of $\Gamma_Q(n)$; that is
$\mbox{Vol}_{m-1}\,  \Gamma_P(n)\geq \mbox{Vol}_{m-1}\, \Gamma_Q(n)$.
Then $m$-dimensional volume of $P$ is not less than
$m$-dimensional volume of $Q$; that is
$\mbox{Vol}_m\, P\geq \mbox{Vol}_m\, Q$.}

Corollary 7.2 can be shortened as follows:.

{\bf 7.3. Corollary.} {\sl Let $P$ and $Q$ be compact convex
polyhedra in $\Bbb R^m$ such that
the volume of every $(m-1)$-dimensional face of $P$ is not less than
the $(m-1)$-dimensional volume of the face of $Q$ parallel to it,
and  the volume of every $(m-1)$-dimensional face of $Q$
is not greater than the
$(m-1)$-dimensional volume of the face of $P$ parallel to it.
Then $m$-dimensional volume of $P$
is not less than the $m$-dimensional volume of $Q$;
that is $\mbox{Vol}_m\, P\geq \mbox{Vol}_m\, Q$.}

{\bf 7.4. Remark.} To emphasize nontriviality of the claims of
7.1--7.3, we note that
under the hypotheses of Corollary 7.3, we cannot say that the
polyhedron $Q$ can be
placed inside the polyhedron $P$ by an appropriate translation.
To demonstrate, take in 3-dimensional space a cube
with edge 10 as $P$
and as $Q$ a rectangular parallelepiped whose faces are parallel
to the faces of $P$
and the edges equal 1, 1 and 50. Then the areas of 2-dimensional
faces of $Q$ are less than
the areas of corresponding parallel faces of $P$, but $Q$ is
``two long'' to be put inside
the cube $P$.
However, according to Corollary 7.3, the volume of the cube
 $P$ is indeed greater than
the volume of the parallelepiped $Q$.

We cannot give explicitly
an analogous example of surfaces with continuous curvature.
Moreover, Blaschke proved that on the plane such an example
is impossible in principle
 \cite{Bl67}.
However, we will show that such an example exists in
$\Bbb R^3$; namely,
there exist infinitely differentiable convex
surfaces
$F, G\subset \Bbb R^3$ such that at the points with parallel
outward normals
the Gaussian curvature of $F$ is always less than or equal to
the Gaussian curvature of $G$
and, nevertheless, $G$ cannot be placed in the convex body
bounded by $F$ by any
translation.

Suppose that the faces of the cube $P$ above are parallel to the
coordinate planes in $\Bbb R^3$.
It is clear that there exists $\alpha$ ($0<\alpha<\pi/2)$
such that for each unit
vector $a\in\Bbb R^3$ there exists a face $p$ of $P$ such
that the outward normal to $p$
makes with $a$ an angle of at most $\alpha$.
Therefore, the area of the projection of $p$ to the plane
orthogonal to $a$ is at least
$100\cos\alpha>0$.
Hence, the area of the projection of the cube onto any plane
in $\Bbb R^3$ is
uniformly separated from 0.

Complicate the situation.
The surface area function $S(P,\cdot)$ of $P$ is concentrated
in the six points
$\pm n_1, \pm n_2, \pm n_3$ on the unit sphere and takes the
value 100 at each of these points.
Given a natural $s$ and an arbitrary unit vector $n$,
denote by
$U_s(n)$ the set of points $x\in\Bbb S^2$ such that the
distance between
$x$ and the endpoint of $n$ is less than $1/s$,
put $N=\{\pm n_1, \pm n_2, \pm n_3\}$ and construct an
infinitely differentiable
function $\nu_s$ on $\Bbb S^2$ such that
(a) $\nu_s$ equals 0 outside the set $\cup_{n\in N} U_s(n)$;
(b) $\nu_s$ is positive in each of $U_s(n)$, $n\in N$;
(c) $\nu_s$ is centrally symmetric, i.e.,
$\nu_s(x)= \nu_s(-x)$ for all $x\in\Bbb S^2$;
(d) the integral of $\nu_s$ with respect to the
standard measure $d\sigma$ on $\Bbb S^2$
over each of $U_s(n)$, $n\in N$, equals 100.

By the Minkowski theorem 4.4, there exists a unique
convex body $P_s$ in $\Bbb R^3$
whose surface area function has density $\nu_s+1/s$:
$$
S(P_s,V)=\int_V (\nu_s+1/s)\, d\sigma, \qquad V\subset\Bbb S^2.
$$

We will check later that each of the bodies $P_s$ is
infinitely differentiable
and the sequence $P_1,\dots, P_s,\dots$ is uniformly bounded
(i.e., lies in a ball of finite radius).
By the Blaschke choice theorem \cite{Bl67},
we can extract a convergent subsequence $P_{s_1},\dots, P_{s_j},\dots$
(convergence is meant with respect to the Hausdorff metric) from the
sequence $P_1,\dots, P_s,\dots$.
Let $P_{s_j}\to P_0$ as $j\to\infty$. Then,
on the one hand, $S(P_{s_j},\cdot)\to S(P_0,\cdot)$,
and on the other hand, by construction,
$S(P_{s_j}, \cdot)\to S(P,\cdot)$. Hence, $S(P_0,\cdot)=S(P,\cdot)$ and,
by uniqueness of a convex body with a given surface area measure,
(see the Minkowski theorem 4.4), $P_0=P$. Therefore,
the sequence of infinitely differentiable surfaces
$P_{s_1},\dots, P_{s_j},\dots$ converges to the cube $P$.

Analogously, we can construct a sequence of infinitely
differentiable surfaces
$Q_{s_1},\dots, Q_{s_j},\dots$ that converges to the
parallelepiped $Q$.
Since the inequality $S(Q,\{n\})<S(P,\{ n\})$ holds for any $n\in N$,
we can smoothen the surface area measure of the
parallelepiped $Q$ so that
the density of the surface area measure of
$Q_{s_j}$ is not greater than the density of the surface
area measure of $P_{s_j}$ (that is than $\nu_{s_j}+1/s_j$).
With this construction, for each $j\in \Bbb N$, the
hypotheses of Theorem~7.1
hold for $P_{s_j}$ and $Q_{s_j}$, however, for all sufficiently large
$j$, there is no translation that shifts
$Q_{s_j}$ inside the body bounded by $P_{s_j}$, because
these surfaces are close to
the boundaries of the polyhedra
$Q$ and $P$, respectively, and the (Euclidean) diameter of
the parallelepiped $Q$ is greater
than the (Euclidean) diameter of the cube $P$.

Let us show now that each of the bodies $P_s$ is infinitely
differentiable.
This follows directly from the following theorem of A.~V.
Pogorelov \cite{Po69}:
{\sl Let $\varkappa$ be a regular $k$ times differentiable
$(k\geq 3)$ positive function
on the unit sphere $\Bbb S^2$.
Let it satisfy the condition
$$
\int_{\Bbb S^2}\frac{x\, d\sigma}{\varkappa (x)}=0,
$$
where $d\sigma$ is the element of area on $\Bbb S^2$ and the
integral is taken over the whole sphere.
Then there exists a regular (at least $k+1$ times differentiable)
surface $F$ with
Gaussian curvature $\varkappa (x)$ at a point with the outward normal $x$.}

Multidimensional variants of this theorem took significant efforts and
generated an extensive literature; for example, see \cite{Bu89}.
These theorems were first  proved by A.~V. Pogorelov \cite{Po71,Po75}.
Later they were re-proved by S.~Y. Cheng and Sh.~T. Yau \cite{CY76};
moreover,
the latter author was awarded the Fields medal for this series of
articles in 1982.
Therefore, sometimes
these multidimensional theorems are referred as
``Pogorelov's theorems for which
Yau has got the Fields medal''.

Finally, let us show that the sequence $P_1,\dots, P_s,\dots$ is
uniformly bounded
(i.e., lies in a Euclidean ball of finite radius).
We will  essentially follow the arguments of
A.~V. Pogorelov \cite{Po69}.
Suppose that the sequence of convex bodies under consideration
is not bounded.
Dropping to a subsequence, we can assume that the (Euclidean) diameter of
$P_s$ is greater than $s$. Then there exists a pair of points
$A_s$ and $B_s$ inside the body
$P_s$ such that the distance between them is greater than $s$.
Project the body $P_s$ onto the plane $\Pi$  orthogonal
to the line
$A_sB_s$. So we obtain a convex figure $\overline{P}_s$.
We have denoted by $\alpha$ ($0<\alpha<\pi/2)$ the number such
that for every unit vector
$a\in\Bbb R^3$ there exists a face $p$ of $P$ such that the
outward normal to $p$
makes with $a$ an angle of at most $\alpha$.
It follows that for every vector $a$ there exists an above domain
$U_s(n)\subset \Bbb S^2$, $n\in N=\{\pm n_1, \pm n_2, \pm n_3\}$
such that an arbitrary vector $x\in U_s(n)$ makes with $a$ an
angle of at most $\alpha -1/s$.
Therefore, the projection of $\widetilde{U_s}(n)\subset\partial P_s$,
which corresponds to the
domain $U_s(n)\subset\Bbb S^2$, onto the plane orthogonal to $a$
has area  separated from 0 uniformly (in $a$ and in $s$):
$$
\mbox{area of the domain projection\, } \widetilde{U_s}(n) =\int\limits_{U_s(n)}\frac{(\nu_s(x)+1/s)\,
d\sigma}{\cos\angle (x,a)}
\geq \frac{100}{\cos (\alpha-1/s)}.
$$

Therefore, for all sufficiently large $s$, the area of the convex
figure $\overline{P}_s$
is not less than some positive $\beta$.
Since every convex figure of diameter $\gamma$ is contained
in some disc
of radius $\gamma$, and so it has area at most $\pi\gamma^2/4$,
the diameter of $\overline{P}_s$ is not less
than $2\sqrt{\beta/\pi}$, that is
not less than some (the same for all sufficiently large $s$)
number $2\sqrt{\beta/\pi}$. Hence, there exist points $\overline{C}_s$ and $\overline{D}_s$
inside the convex figure $\overline{P}_s$ such that the
Euclidean distance between them
is at least $\sqrt{\beta/\pi}$. Let $C_s$ and $D_s$ be
points of a convex body bounded by
$P_s$ which are projected to the points
$\overline{C}_s$ and $\overline{D}_s$, respectively, on the plane $\Pi$.
Construct a plane $\widetilde{\Pi}$ parallel to the lines $
A_sB_s$ and $C_sD_s$.
Project the surface $P_s$ to the plane $\widetilde{\Pi}$.
This projection contains
a quadrilateral whose vertices are the projections of
$A_s$, $B_s$, $C_s$ and $D_s$. The area of this quadrilateral
is obviously not less
than the product of the lengths of the segments $A_sB_s$ and
$C_sD_s$, that is not less than
$s\sqrt{\beta/\pi}$. Therefore, taking $s$ to be sufficiently
large, we can make
the area of the projection of $P_s$ onto
$\widetilde{\Pi}$ however large.
Hence, taking
$s$ to be sufficiently large, we can make the area of
the surface $P_s$ however larger. The latter case is
impossible since the convex surfaces
$P_s$ converge to the surface $P$ of finite area.
This contradiction completes the proof
of the fact that there exist smooth surfaces in
$\Bbb R^3$ satisfying the hypotheses of Theorem~7.1 which
cannot be shifted inside each
other by any translation. A multidimensional example can be
constructed analogously.

{\bf 7.5. Theorem.} {\sl For every two convex bodies $K$ and
$L$ in $\Bbb R^m$ $(m\geq 2)$,
the volume of their Blaschke sum does not exceed the volume
of their Minkowski sum:
$\mbox{Vol\, }(K\# L)\leq\mbox{Vol\, }(K+L)$.}

{\bf Proof.} First of all, note that it is sufficient to prove
the theorem for polyhedra. The general case is then
obtained by passage to a limit.

For polyhedra, the proof is by induction by the dimension of
the space $\Bbb R^m$.
As it was noticed in Section~5.1, the Blaschke sum of convex
polygons coincides with
their Minkowski sum on the plane: $K\# L=K+L$.
Therefore, for $m=2$ the claim of Theorem 7.5 is trivial.

If $m=3$ then it is known \cite{Ly56} that the
Minkowski sum of two convex
polyhedra is a convex polyhedron whose faces are obtained as
the
Minkowski sum of
(a)the  faces of the  polyhedra under addition;
(b) a face of one of these
polyhedra and
an edge or a vertex of the other;
(c) nonparallel edges of the polyhedra. Here, all mentioned
faces, edges and vertices
lie  in the supporting planes with parallel outward normals.
The cases (a) and (b)
are especially important for us.

Indeed, if the summands of the Minkowski sum have faces with
the same normal, then
the sum has a face with the same normal and area greater than
or equal to the sum
of areas of the summands; since, by (a), the parallel faces
$\Gamma_1$ and $\Gamma_2$ must be added in the sense of
Minkowski and, by the Brunn--Minkowski
inequality (3),
$$
\mbox{Vol\, }(\Gamma_1+\Gamma_2)^{1/2}\geq \mbox{Vol\, }\Gamma_1^{1/2}+\mbox{Vol\, }\Gamma_2^{1/2}.
$$
It follows that
$$
\biggl[\frac{\mbox{Vol\, }(\Gamma_1+\Gamma_2)}{\mbox{Vol\, }\Gamma_1}\biggr]^{1/2}\geq
1+\biggl[\frac{\mbox{Vol\, }\Gamma_2}{\mbox{Vol\, }\Gamma_1}\biggr]^{1/2},
$$
and so
$$
\frac{\mbox{Vol\, }(\Gamma_1+\Gamma_2)}{\mbox{Vol\, }\Gamma_1}\geq
\biggl\{1+\biggl[\frac{\mbox{Vol\, }\Gamma_2}{\mbox{Vol\, }\Gamma_1}\biggr]^{1/2}\biggr\}^2\geq
1+\frac{\mbox{Vol\, }\Gamma_2}{\mbox{Vol\, }\Gamma_1}
$$
or
$$
\mbox{Vol\, }(\Gamma_1+\Gamma_2)\geq\mbox{Vol\, }\Gamma_1+\mbox{Vol\, }\Gamma_2.
$$
Therefore, if we add, in the sense of Minkowski, the faces with
the same outward normals,
then we obtain a face with the same normal and, moreover, the
area of this face is
not less than the sum of the areas of the summands. That is
in this case
the face of $K+L$ has area  not less than
the corresponding face
of $K\# L$ parallel to it.

Considering the cases (b) and (c) analogously, we see that
the area of each face of
the Minkowski sum is not less than the area of a face of the
Blaschke sum $K\# L$
parallel to it. Hence, Theorem~7.5 for polyhedra follows from
Theorem~7.1 (and its Corollaries~7.2 or~7.3).

In the general case, the induction step  is analogous to the
above step from $m=2$ to $m=3$.

{\bf 8. On exponents in the Brunn--Minkowski and
Kneser--S\"uss inequalities.}
In the proof of Theorem~7.5 we have given, in fact, an
``arithmetical'' reason for the fact
that the concavity of the square root of the Minkowski sum
area for plane figures
implies the concavity of the Minkowski sum area itself.
This suggests that, firstly, exponents in the
Brunn--Minkowski and
Kneser--S\"uss inequalities can be changed in some range so
that the inequalities
keep, and, secondly, that there must be an ``unimprovable'' exponent.
As far as we know, the question was not posed in this form
before and we
think that it is worth to give here some (although very simple)
arguments for that.

{\bf 8.1. Theorem.} {\sl For every $a\geq 1$ and all convex bodies $K$
and $L$
in $\Bbb R^m$ $(m\geq 2)$ the following inequalities hold:
$$
(\mbox{Vol\, }(K+L))^{a/m}\geq (\mbox{Vol\, }K)^{a/m}+(\mbox{Vol\, }L)^{a/m} \eqno(4)
$$
$$
(\mbox{Vol\, }(K\#L))^{a-a/m}\geq (\mbox{Vol\, }K)^{a-a/m}+(\mbox{Vol\, }L)^{a-a/m}. \eqno(5)
$$
Moreover, for  $0<a<1$ there exist convex bodies $K$ and $L$ in
$\Bbb R^m$ $(m\geq 2)$ such that the inequalities
$(4)$ and $(5)$ fail.}

In other words, wee can say that, increasing exponent
in the Brunn--Minkowski or Kneser--S\"uss
inequalities, we coarsen the inequalities but they still hold.
The original Brunn--Minkowski and Kneser--S\"uss inequalities are
optimal in the sense that
the exponents cannot be decreased.

For the proof of Theorem~8.1 we  need the following absolutely
elementary lemma whose proof we leave to the reader.

{\bf 8.2. Lemma.} {\sl  For all $a\geq 1$ and $x>0$ the following
inequality holds:
$$
(1+x)^a\geq 1+x^a.\eqno(6)
$$
Moreover, for all $0<a<1$ there exists $x>0$ such that $(6)$ fails.}

{\bf Proof of Theorem~8.1} is given only for the Brunn--Minkowski
inequality,
since for the Kneser--S\"uss inequality the same arguments
can be applied.

For  $a\geq 1$, by the Brunn--Minkowski inequality and
the inequality~(6), we have
$$
\biggl[\frac{\mbox{Vol\, }(K+L)}{\mbox{Vol\, }K}\biggr]^{a/m}\geq
\biggl\{1+\biggl[\frac{\mbox{Vol\, }L}{\mbox{Vol\, }K}\biggr]^{1/m}\biggr\}^a\geq
1+\biggl[\frac{\mbox{Vol\, }L}{\mbox{Vol\, }K}\biggr]^{a/m},
$$
that proves (4) for $a\geq 1$.

Now, let  $0<a<1$. Take a ball of radius $r$ as a body $K$ and
a ball of radius $R$ as $L$.
Then $K+L$ is a ball of radius $r+R$ and we have
$\mbox{Vol\, }K=\omega_m r^m$, $\mbox{Vol\, }L=\omega_m R^m$
and $\mbox{Vol\, }(K+L)=\omega_m (r+R)^m$,
where $\omega_m$ is the volume of the unit ball in $\Bbb R^m$.
Then (4) is equivalent to
$$
\bigl[\omega_m (r+R)^m\bigr]^{a/m}\geq  \bigl[\omega_m r^m\bigr]^{a/m} +  \bigl[\omega_m R^m\bigr]^{a/m},
$$
or
$$
(r+R)^a\geq  r^a + R^a,\qquad\mbox{£³£}\qquad \biggl(1+\frac{r}{R}\biggr)^a\geq 1+\biggl[\frac{r}{R}\biggr]^a.
$$
However, by Lemma~8.2 for $0<a<1$, the latter inequality
certainly fails for some $x=r/R$.

{\bf 9. Visualization of the Blaschke sum of convex polyhedra.}
On the base of the package OpenGeometry \cite{GS99},
the authors developed a
computer program
for visualization of the Blaschke sum of polyhedra in 3-dimensional
Euclidean space.
The main difficulty is to construct a polyhedron
given by a set of outward normals
and areas of faces. Let us describe an algorithm for that.

Suppose that we plan to construct a polyhedron $P$ such that
$n_1,\dots,n_k$ are unit vectors of outward normals to its faces and
$F_1,\dots,F_k$ are areas of faces.

We start with a polyhedron $P_0$ that circumscribes the
unit sphere with the center
the origin and has outward normals $n_1,\dots,n_k$, then we deform
it by modifying
the support numbers so that the vector
$n_1,\dots,n_k$ remain unchanged.

{\bf Step 1.} Let $P_0$ be a polyhedron with outward normals
$n_1,\dots,n_k$ and support numbers $h_j=1$ for all $j=1,\dots,k$.
It is clear that such a polyhedron always exists. Denote by
$F_1^0,\dots,F_k^0$ areas of faces of $P_0$.  Then
$\sum\limits_{j=1}^k F_j^0 n_j=0$.

{\bf Step 2.} For every $t\in[0,1]$, consider a polyhedron $P_t$ which is
the Blaschke sum $(1-t)P_0\# tP$, i.e., a polyhedron with outward normals
$n_1,\dots,n_k$ and the area of the $j$-th face is given by the formula
$F_j^t=(1-t)F_j^0+tF_j$. Clearly, for every $t\in[0,1]$, we have
$\sum\limits_{j=1}^k F_j^t n_j=0$ and so $P_t$ exists.

{\bf Step 3:} Computing the support numbers $h_1^t,\dots,h_k^t$
of the planes
$\alpha_1^t,\dots,\alpha_k^t$ that bound $P_t$.
Partition the interval $[0,1]$ into equal intervals of length $\Delta t$.
Suppose that the support numbers $h_1^t,\dots,h_k^t$ of $P_t$
 have been computed
for the left endpoint of an interval of length $\Delta t$.
The support numbers $h_1^{t+\Delta t},\dots,h_k^{t+\Delta t}$ of
the polyhedron $P_{t+\Delta t}$
that corresponds to the right endpoint of this interval will be found
from a linear
system of algebraic equations. Since for $t=0$ all support numbers
are known
(they are equal to 1), the problem will be solved at $t=1$.

Let us describe the construction for a required system of equations.

Each polyhedron $P_t$ is determined by its support numbers
$h_1^t,\dots,h_k^t$. Therefore, the set of all polyhedra $P_t$ can be
presented as a subset of $k$-dimensional space ${\Bbb R}^k$
with coordinates $h_1^t,\dots,h_k^t$.
This subset is open, because the faces do not disappear under small
displacement.

Combine in a single class all polyhedra $P_t$ that are equal and parallel
to each other.
Since each translation is determined by three components,  such
a class is
determined by $k-3$ variables.
The set of these classes forms a $(k-3)$-dimensional manifold~${\bf A}_t$.

Consider the set ${\bf B}_t$ of all $k$-tuples of numbers
$F_1^t,\dots,F_k^t$ such that $F_j^t>0$ for all $j=1,\dots,k$ and
$\sum_{j=1}^k F_j^t n_j=0$. The vectors $n_1,\dots,n_k$ are assumed to be
fixed.
The vector identity $\sum_{j=1}^k F_j^t n_j=0$
is equivalent to three scalar identities and so it determines a
$(k-3)$-dimensional plane in ${\Bbb R}^k$. The inequalities $F_j^t>0$
determine in
this plane an open convex set which is exactly ${\bf B}_t$.

Thus, we have a natural one to one map
$F:{\bf A}_t\to {\bf B}_t$ (see \cite{AD50}).

Let
$\Delta F_j=F_j^{t+\Delta t}-F_j^t$ and $\Delta h_j=h_j^{t+\Delta t}-h_j^t$.
Then
$$
\Delta F_j=\frac{\partial F_j}{\partial h_1} \Delta h_1+\dots+
\frac{\partial F_j}{\partial h_k} \Delta h_k, \quad j=1,\dots,k. \eqno(7)
$$
Note that $\Delta F_j=\Delta t(-F_j^0+F_j)$, i.e., $\Delta F_j$
is a constant for a given
$\Delta t$.

We are left with computing $\frac{\partial F_j}{\partial h_i}$.
Let $j$-th face of $P_t$ is cut by planes
$\alpha_{q_1},\dots,\alpha_{q_{m_j}}$. Denote by $\ell_{jp}$
the length of an edge lying on the line $\alpha_j\cap\alpha_p$.
It is not difficult to show that the modification of the support
number $h_j$ implies the modification of the area
of the $j$-th face to the value
$$
\frac{\partial F_j}{\partial h_j}=\sum\limits_{p\in\{q_1,\dots,q_{m_j}\}}\ell_{jp}
{\rm ctg}\angle(n_j,n_p),
$$
and the modification of the area of the $i$-th face ($i\not=j$)
to the value
$$
\frac{\partial F_i}{\partial h_j}=\frac{\ell_{ji}}{\sin\angle(n_j,n_i)}.
$$
Hence, the system (7) is a system of linear equations of rank
$k-3$ in $k$ unknowns $\Delta h_1,\dots,\Delta h_k$.
We take three of the unknowns arbitrarily (for example,  equal to~0)
and find the other $k-3$ unknowns from~(7).
Thus, we obtain $\Delta h_1,\dots,\Delta h_k$, and so
$h_j^{t+\Delta t}=\Delta h_j+h_j^t$.

{\bf Step 4.} We now  construct $P_t$ with the support
numbers $h_1^t,\dots,h_k^t$.
In order to do this, given every pair of planes
$\alpha_i^t$ and $\alpha_j^t$, we find a set of points
$\{\alpha_i^t\cap\alpha_j^t\cap\alpha_p^t, p=1,\dots,k\}$ and then
we pick up from this set those points that are vertices of the polyhedron.
Since the polyhedron is convex, there are at most two of these points.
We add these
points to the list of vertices $\{V\}$ and, if the points were two,
put the elements $\ell_{ij}$ and $\ell_{ji}$ to be equal to the
distance between
these points.
In all other cases we put $\ell_{ij}=\ell_{ji}=0$.

When the list $\{V\}$ is composed,  we, for each plane
$\alpha_i^t$, pick up from $\{V\}$ those points that lie in $\alpha_i^t$.
So we determine the  $i$-th face of $P_t$.

\medskip

\begin{figure}
\begin{minipage}[b]{0.30\linewidth}
\begin{verbatim}
10
1 1 0 2.0412414523
0 1 1 2.0412414523
1 0 1 2.0412414523
-1 1 1 5
1 -1 1 5
1 1 -1 5
-1 -1 1 5
-1 1 -1 5
1 -1 -1 5
-1 -1 -1 5
\end{verbatim}
\centerline{a)}
\end{minipage}\hfill
\begin{minipage}[b]{0.30\linewidth}
\centering
\includegraphics[width=3.2 cm]{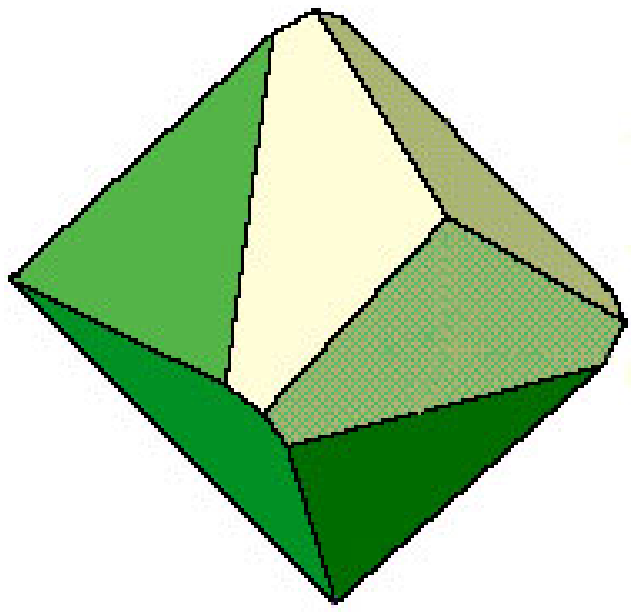}\\
\bigskip
b)
\end{minipage}\hfill
\begin{minipage}[b]{0.30\linewidth}
\centering
\includegraphics[width=3.2 cm]{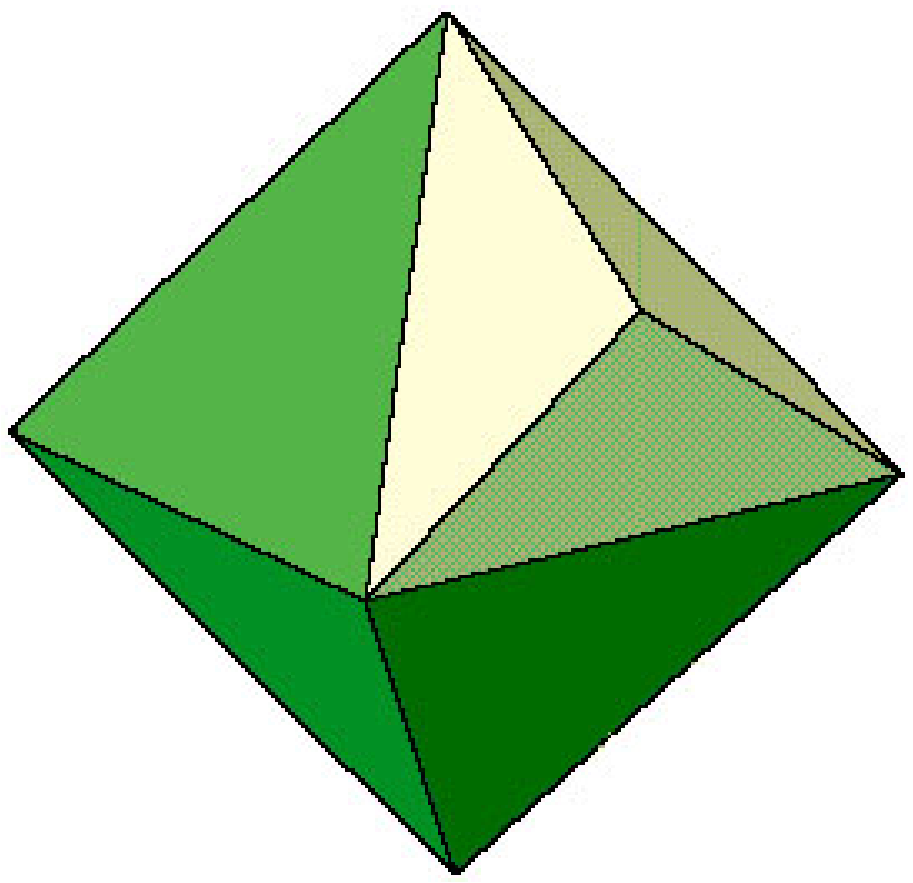}\\
\bigskip
c)
\end{minipage}\hfill
\caption{Construction of the Gr\"unbaum polyhedron}\label{grunbaum}
\end{figure}

{\bf 9.1. Example.} In the process of deformation, $P_t$ can
change its combinatorial structure.
In order to be convinced that the proposed algorithm works correctly
with disappearing and appearing edges, consider a polyhedron $P$
obtained from an octahedron by attaching a tetrahedron with
a face congruent to a
face of the octahedron (see Figure~\ref{grunbaum}~c)).

The polyhedron $P$ is input to the programme as shown
in Figure~\ref{grunbaum}~a).
Here, the first line gives the number of faces and each
subsequent line
shows coordinates of an outward (not necessarily unit) normal
to a face
and the area of this face.

In \cite{Gr63} and \cite[p. 286]{Gr03}, it was shown that a
polyhedron of such combinatorial
type cannot circumscribe a sphere. A polyhedron
$P_0$ that circumscribes a sphere with the same outward normals as $P$
is shown in Figure~\ref{grunbaum} b). The difference in the
combinatorial structures of
$P$ and $P_0$ is obvious.

{\bf 9.2. Remark}. Example 9.1 shows, in particular, that it is
impossible to determine
a combinatorial structure of a polyhedron if only outward normals
are known.
(That is it is impossible to determine which faces of the polyhedron
have a common edge.)

{\bf 9.3. Example:} The Blaschke sum of a cube and icosahedron.
In Figure~\ref{cube_ico}~c), the Blaschke sum of a cube whose
faces are parallel
o the coordinate planes and face areas equal 2 and an icosahedron
given as shown in
Figure~\ref{ico_dodeca} is shown.

\begin{figure}
\begin{minipage}[t]{.6\linewidth}
Icosahedron:
\begin{verbatim}
20
1 1 1 5
-1 1 1 5
1 -1 1 5
1 1 -1 5
-1 -1 1 5
-1 1 -1 5
1 -1 -1 5
-1 -1 -1 5
0 0.6180339887 1.618033989 5
0 -0.6180339887 1.618033989 5
0 0.6180339887 -1.618033989 5
0 -0.6180339887 -1.618033989 5
0.6180339887 1.618033989 0 5
-0.6180339887 1.618033989 0 5
0.6180339887 -1.618033989 0 5
-0.6180339887 -1.618033989 0 5
1.618033989 0 0.6180339887 5
-1.618033989 0 0.6180339887 5
1.618033989 0 -0.6180339887 5
-1.618033989 0 -0.6180339887 5
\end{verbatim}
\end{minipage}\hfill
\begin{minipage}[t]{.4\linewidth}
Dodecahedron:
\begin{verbatim}
12
0 1.618033989 1 3
0 1.618033989 -1 3
0 -1.618033989 1 3
0 -1.618033989 -1 3
1 0 1.618033989 3
-1 0 1.618033989 3
1 0 -1.618033989 3
-1 0 -1.618033989 3
1.618033989 1 0 3
-1.618033989 1 0 3
1.618033989 -1 0 3
-1.618033989 -1 0 3
\end{verbatim}
\end{minipage}\hfill
\caption{Input data for an icosahedron and dodecahedron}\label{ico_dodeca}
\end{figure}

\begin{figure}
\centering
\begin{tabular}{ccccc}
\raisebox{0.5cm}{\includegraphics[width=1.6 cm]{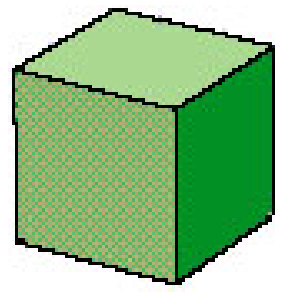}}&
\raisebox{1.3cm}{$\#$} &
\includegraphics[width=3 cm]{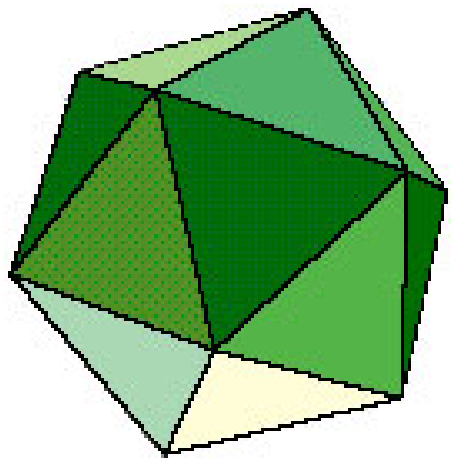} &
\raisebox{1.3cm}{$=$} &
\includegraphics[width=3 cm]{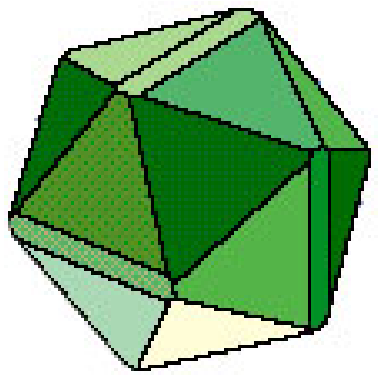}\\
a) & & b) & & c)
\end{tabular}
\caption{The Blaschke sum of a cube and icosahedron}\label{cube_ico}
\end{figure}

The result enables us to conjecture that an icosahedron has parallel
edges, which is,
of course, true, and that there is a simple way to construct
an icosahedron.

\begin{figure}
\centering
\begin{tabular}{ccc}
\includegraphics[width=2.5 cm]{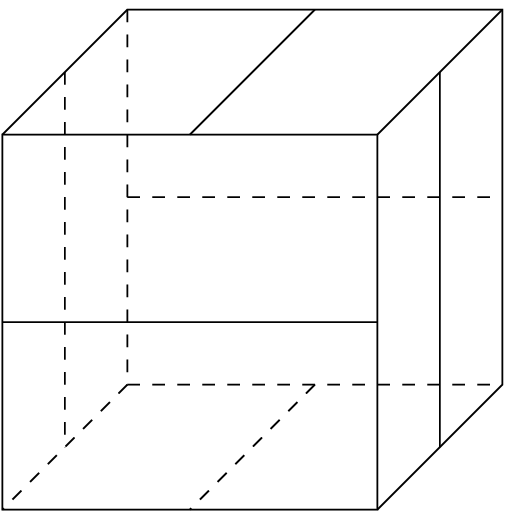}\quad & \quad
\includegraphics[width=2.5 cm]{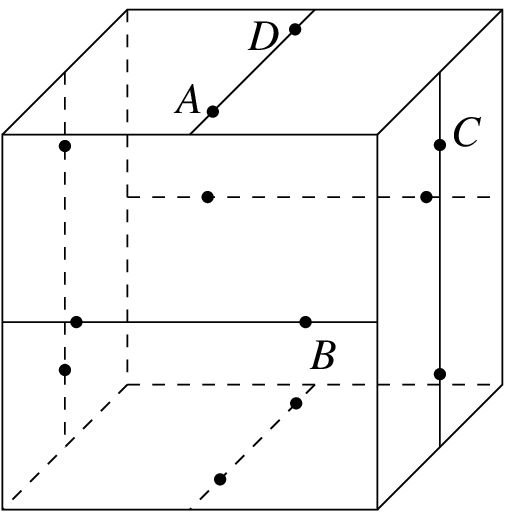}\quad & \quad
\includegraphics[width=2.5 cm]{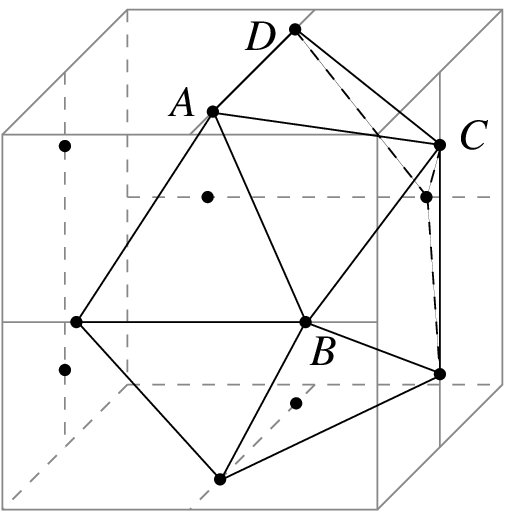}\\
a) \quad & \quad b) \quad & \quad c)
\end{tabular}
\caption{Construction of an icosahedron}\label{ico_constr}
\end{figure}

Take a cube with an edge 2 and construct segments with the endpoints
at the edges
as shown in Figure~\ref{ico_constr}~a). On each segment, mark two
points at the distance
$t=(3-\sqrt{5})/2$ from those edges of the cube to which this
segment is orthogonal
(Figure~\ref{ico_constr}~b)).
Cutting parts off the cube by planes that pass through triples of
points according
to the combinatorial structure of an icosahedron (for example,
by planes through
$A$, $B$ and $C$, and through $A$, $C$ and $D$), we get an icosahedron
(Figure~\ref{ico_constr}~c)).

\begin{figure}
\centering
\begin{tabular}{ccccc}
\raisebox{0.8cm}{\includegraphics[width=2 cm]{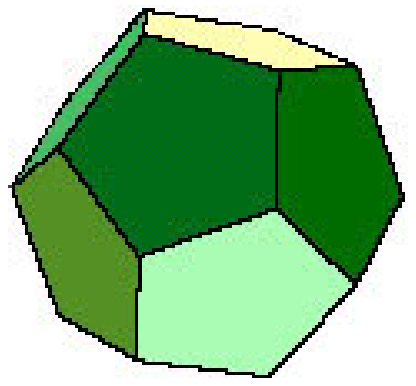}}&
\raisebox{1.7cm}{$\#$}&
\raisebox{0.2cm}{\includegraphics[width=3 cm]{icosahedron.eps}} &
\raisebox{1.7cm}{$=$} &
\includegraphics[width=3.3 cm]{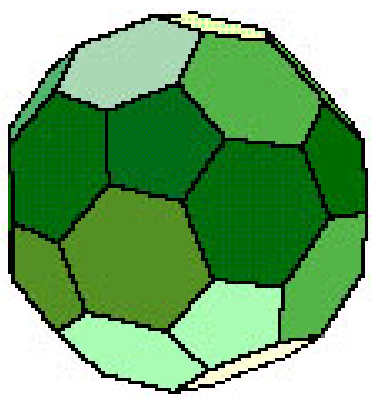}\\
a) & & b) & & c)
\end{tabular}
\caption{The Blaschke sum of a dodecahedron and icosahedron}\label{football}
\end{figure}

{\bf 9.4. Example:} The Blaschke sum of a dodecahedron and
icosahedron.
In Figure~\ref{football}~c), the Blaschke sum of a dodecahedron
and icosahedron
given as in Figure~\ref{ico_dodeca} is shown. Shortly
speaking, the Blaschke sum of
a dodecahedron and icosahedron is a football
(which is well known to be sewed from flat leather penta- and
hexagons so that
each pentagon is sewed only with hexagons and pentagons and
hexagons are attached to a hexagon
consequently).
Although, the polyhedron shown in Figure~\ref{football}~c) has
a quite scientific name:
``truncated icosahedron'' (see \cite{We74}).

{\bf 9.5. Remark.}
Our program for visualization of a polyhedron with given areas of
faces and outward normals
is not unique. Something like that was done in \cite{Zo03a,Zo03b}
for the purposes of
computer graphics. Unfortunately, the papers \cite{Zo03a,Zo03b} are
not available for us.

{\bf 9.6. Remark.}
Our program for visualization of the Blaschke sum enables us to find,
as a byproduct,
the numerical values of various geometrical parameters of the
constructed
polyhedra, for example,
the volume, areas of faces, the integral mean curvature etc.
This enables us to master some geometric intuition and reject
some conjectures that arose during the work.

\bigskip

Victor Alexandrov

Natalia Kopteva

Semen S. Kutateladze

\medskip

Sobolev Institute of Mathematics

acad. Koptyug ave. 4,

630090 Novosibirsk, Russia

E-mail address: \{alex, natasha, sskut\}@math.nsc.ru

\end{document}